\documentclass[12pt,twoside]{amsart}
\usepackage{amssymb}

\nonstopmode

\numberwithin{equation}{section}
\font\tengothic=eufm10 scaled\magstep 1
\font\sevengothic=eufm7 scaled\magstep 1
\newfam\gothicfam
          \textfont\gothicfam=\tengothic
          \scriptfont\gothicfam=\sevengothic

\newcommand {\PP}{\mathbb{P}}
\newcommand{\cI}{\mathcal{I}}

\newcommand{\cL}{\mathcal{L}}
\newcommand{\cN}{\mathcal{N}}
\newcommand{\cE}{\mathcal{E}}
\newcommand{\cO}{\mathcal{O}}

\DeclareMathOperator{\pnt}{\raise 0.5mm \hbox{\large\bf.}}

\newcommand{\s}{\; | \;}
\newtheorem{theorem}{Theorem}[section]
\newtheorem{lemma}[theorem]{Lemma}
\newtheorem{proposition}[theorem]{Proposition}
\newtheorem{corollary}[theorem]{Corollary}

\theoremstyle{definition}
\newtheorem{definition}[theorem]{Definition} 
\newtheorem{remark}[theorem]{Remark}
\newtheorem{example}[theorem]{Example}

\begin{document}


\title[Liaison addition and the structure of a Gorenstein liaison class]
{Liaison addition and the structure of a Gorenstein liaison class}

\author{Robin Hartshorne}
\author{Juan Migliore}
\author{Uwe Nagel}

\address{Department of Mathematics,
         Evans Hall,
         University of California,
         Berkeley, CA 94720-3840,
         USA}
\email{robin@math.berkeley.edu}

\address{Department of Mathematics,
         University of Notre Dame,
         Notre Dame, IN 46556,
         USA}
\email{migliore.1@nd.edu}
\address{Department of Mathematics,
         University of Kentucky,
         715 Patterson Office Tower,
         Lexington, KY 40506-0027,
         USA}
\curraddr{Institute for Mathematics \& its Applications,
          University of Minnesota,
          Minneapolis, MN 55455-0436,
          USA}
\email{uwenagel@ms.uky.edu}

\begin{abstract}
We study the concept of liaison addition for codimension two subschemes of an arithmetically Gorenstein projective scheme. We show how it relates to liaison and biliaison classes of subschemes and use it to investigate the structure of Gorenstein liaison equivalence classes, extending the known theory for complete intersection liaison of codimension two subschemes. In particular, we show that on the non-singular quadric threefold in projective 4-space, every non-licci ACM curve can be obtained from a single line by successive liaison additions with lines and CI-biliaisons.
\end{abstract}

\maketitle

\section{Introduction}

The Lazarsfeld-Rao property refers to a structure common to even
liaison classes in codimension two under complete intersection
liaison (cf.\ \cite{BBM}, \cite{MP}, \cite{BM4},
\cite{Bolondi-Migliore-manuscripta}, \cite{N-gorliaison},
\cite{RTLRP}). The goal of our work is to discover if there is an
analogue of this property for Gorenstein liaison.  In other words,
we seek to find some structure for an even Gorenstein liaison class.

At this point, more precisely, the Lazarsfeld-Rao theorem for
complete intersection liaison is known only for subschemes of
codimension two in an ambient projective scheme $X$, but $X$ can be
taken quite generally: the most general result to date is for $X$ an
integral projective scheme satisfying the condition $S_3$ of Serre
and $H^1_*({\mathcal O}_X) = 0$ (\cite{RTLRP}).  In this case it
says:

\begin{enumerate}
\item If $C$ is a codimension two subscheme of $X$ (equidimensional
without embedded components) that is not of minimal degree in its
CI-biliaison equivalence class, then it admits a strictly descending
biliaison.

\item Any two subschemes $C,C'$ of minimal degree in the same biliaison equivalence
class can be joined by a sequence of elementary biliaisons of height 0.
\end{enumerate}

Since nothing has been proved for subschemes of codimension three or
higher in any ambient scheme, we will also stick to codimension two.
We refer to \cite{CH} and \cite{CDH} for definitions and basic results on
CI-liaison, Gorenstein liaison, and Gorenstein biliaison on a projective
scheme X. Our basic assumption throughout this paper is that $X \in \PP^N$ is a
normal arithmetically Gorenstein scheme, and that we deal with closed
subschemes $C$ that are equidimensional of codimension 2 without embedded points.
Recall that a coherent sheaf $\mathcal N$ on $X$ is called {\em extraverti} \cite[Definition 2.9]{RTLRP} if
$H^1_*({\mathcal N}^\vee) = 0$ and ${\mathcal E}xt^1({\mathcal
N},{\mathcal O}_X) = 0$. A sheaf  ${\mathcal F}$ is {\em dissoci\'e} if ${\mathcal F} = \bigoplus
{\mathcal O}_X(a_i)$ for some integers $r, a_1,\ldots,a_r$. An {\it ${\mathcal N}$-type resolution} of $C$  \cite[Definition 2.4]{CDH} is an
exact sequence
\[
 0 \to {\mathcal L} \to {\mathcal N} \to I_C \to 0
\]
with ${\mathcal L}$ dissoci\'e and ${\mathcal N}$ extraverti. In addition if $C$ is locally
Cohen-Macaulay, then ${\mathcal N}$ is a locally Cohen-Macaulay sheaf on  $X$.
There are two
(inequivalent) kinds of Gorenstein liaison that will concern us, so
we phrase two questions.

\medskip

\noindent {\bf Question 1}: What is the structure of a Gorenstein
biliaison equivalence class of codimension two subschemes of $X$
(see \cite{CH} for Gorenstein biliaison)?

\medskip

\noindent {\bf Question 2}: What is the structure of an even
Gorenstein liaison class of codimension two subschemes of $X$?

\medskip

Since any CI-biliaison is also an even Gorenstein liaison and a
Gorenstein biliaison, the first observation we can make is that each
Gorenstein biliaison or even Gorenstein liaison class is a disjoint
union of CI-biliaison classes, and within each of these, the
``classical" Lazarsfeld-Rao property holds.  So our problem is
rather, how to get from one CI-biliaison class to another within
the Gorenstein biliaison or even Gorenstein liaison class.  And,
while we are at it, do our constructions yield minimal elements of a
CI-biliaison class if we start from a minimal element?

Interesting special cases are the following classes, all  contained
within the set of arithmetically Cohen-Macaulay subschemes:
\[
\hbox{\{licci\}} \subseteq \hbox{\{gobilicci\}} \subseteq \hbox{\{glicci\}}.
\]
Here we follow the commonly used acronyms for the liaison  class
(resp.\ Gorenstein biliaison class; resp.\ Gorenstein liaison class)
of a complete intersection.

Another general observation is that in the case of Gorenstein
biliaison, we have elementary Gorenstein biliaisons, which could
play the role of elementary biliaisons in the traditional
Lazarsfeld-Rao property.  But in the case of even Gorenstein
liaison, it may happen that the only Gorenstein biliaisons are
CI-biliaisons, so we need some other operation to move from one
class to another. We investigate liaison addition as a possibility
in this case, and have some success.  In particular, we give some
additional justification to the name itself.

``Liaison Addition" was introduced by Schwartau in his thesis,
where on page 1 he says: ``Does there exist a geometric addition of
curves in ${\mathbb P}^3$ corresponding to the direct sum of their
liaison invariants?  This is the \underline{liaison addition}
problem. \dots We find that not only is there a way to add curves in
${\mathbb P}^3$, but that an explicit procedure is possible: that
is, equations for the added curve may be written down from the
equations of the curves being added.  The addition procedure \dots
admits a purely intrinsic formulation reminiscent of liaison
itself."

Since Schwartau considered the question only in ${\mathbb P}^3$,
where Rao had shown that liaison reduces to a question about the
``liaison invariants" (subsequently dubbed ``Rao modules"), the name
made perfect sense in his setting.  Subsequently, liaison theory has
exploded in the direction of Gorenstein liaison, thanks largely to \cite{KMMNP}
(see \cite{book} for
an extensive bibliography, albeit now quite outdated).  Liaison
addition has also been generalized substantially (see \cite{BM4}, \cite{Bolondi-Migliore-manuscripta}, \cite{GM4}, \cite{book}); a treatment in the generality needed here can be found, for example, in \cite{N-gorliaison}.
The name has continued to make sense in the context of complete
intersection liaison in codimension two (cf.\
\cite{Bolondi-Migliore-manuscripta}, \cite{N-gorliaison}).  However,
until now there has been no connection made between liaison addition
and the more general notion of Gorenstein liaison.


In section 2 of this paper we recall the construction and first properties
of liaison addition.  Our first main results of this paper are
contained in section 3, where we prove the following about the liaison addition $C$
 of given codimension two
subschemes $C_1$ and $C_2$ with respect to forms $F_1 \in I_{C_2}$
and $F_2 \in I_{C_1}$.

\begin{itemize}
\item If $C_2$ is gobilicci then $C_1$ and $C$ are in the
same   Gorenstein biliaison class on~$X$.

\item If $C_2$ is glicci then $C_1$ and $C$ are in the
same even Gorenstein liaison class on $X$.
\end{itemize}

It follows from  Rao's
theorem and the preparatory results on liaison addition that if
$C_2$ is licci then $C_1$ and $C$ are in the same CI even liaison
class.  However, it is only an existence result about a sequence of
links.  In section 4 we make this more precise by showing that there is a very concrete sequence of
links, preserving the liaison addition structure all along the way
from $C$ to $C_1$.

In section 5 we begin the study of a Lazarsfeld-Rao-type  structure
theorem for a Gorenstein liaison class.  Two specific situations
that we will examine in some detail are the case of curves on the
non-singular quadric hypersurface in ${\mathbb P}^4$ and on the
singular quadric hypersurface with a single double point in
${\mathbb P}^4$.  We hope that these cases will illustrate the type
of phenomena one finds, and may suggest what kind of results one
could hope for in more general situations.  The special feature of
these two examples is that in each case, the Rao module of a curve
characterizes the even Gorenstein liaison (resp.\ G-biliaison)
equivalence class of a curve, in analogy to the traditional Rao
theorem, where the Rao module characterizes the CI-biliaison class
of a curve in ${\mathbb P}^3$.  (See \cite[Theorem 6.2]{CDH} for
the first case and \cite[Theorem 6.2]{CH} for the second.)  Note
that it is an open question whether the Rao module characterizes an
even Gorenstein liaison class of curves in ${\mathbb P}^4$.


\section{Liaison addition}

Throughout this note we denote by $R = k[x_0,\ldots,x_n]$ the homogenous coordinate ring of $\PP^n$ where $k$ is any infinite field. We begin with the definition.

\begin{definition} \label{der-liai-add}
Let $C_1, C_2$ be codimension 2 subschemes of $X  \subset \PP^n$.  Let $F_1
\in I_{C_1, X}$ and $F_2 \in I_{C_2, X}$ be homogenous elements of
degree $f_1$ and $f_2$, respectively, such that $\{F_1, F_2\}$ is an
$S$-regular sequence, where $S = R/I_X$. Then the subscheme $C \subset X$ defined by the ideal
\[
I_{C, X} := F_2 \cdot I_{C_1, X} + F_1 \cdot I_{C_2, X}
\]
is called the {\em liaison addition of $C_1$ and $C_2$ with respect to $F_1$ and $F_2$}.
\end{definition}

We record some of its properties:

\begin{lemma} \label{lem-liai-add}
Assume $X$ is an arithmetically Gorenstein scheme and
that $C_1$ and $C_2$ are codimension 2 equidimensional subschemes.
Then the ideal $J = F_2 \cdot I_{C_1, X} + F_1 \cdot I_{C_2, X}$ is a
saturated ideal in the coordinate  ring $S = R/I_X$ of $X$.  Thus,
it is the homogeneous ideal of a subscheme $C \subset X$ which has
the following properties:
\begin{itemize}
\item[(a)] The Hilbert function of $C$ is for all integers $j$:
$$
h_C (j) = h_{C_1}(j-f_2) + h_{C_2} (j - f_1) + h_Y (j)
$$
where $Y$ is the complete intersection defined by $(F_1, F_2)$.
\item[(b)] Let
$$
0 \to \cL_i \to \cN_i \to \cI_{C_i, X} \to 0
$$
be an $\mathcal N$-type resolution of $C_i$, $i = 1, 2$, on $X$.
Then  $C$ has the following $\mathcal N$-type resolution on $X$
$$
0 \to \cO_X (-f_1 - f_2) \oplus \cL_1 (-f_2) \oplus \cL_2 (-f_1) \to
\cN_1 (-f_2) \oplus \cN_2 (-f_1) \to \cI_{C, X} \to 0.
$$
\item[(c)] If $\dim X = 3$ with $d = \deg X$ and $\omega_X = {\mathcal O}_X(e)$, and
$C_1, C_2$ are curves of degrees $d_1$ and $d_2$ and arithmetic genera $g_1, g_2$,
respectively, then the  degree of $C$ is
$$
\deg C = d_1 + d_2 + d f_1 f_2
$$
and its arithmetic genus is
$$
g_C = g_1 + g_2 - 1 + d_1 f_2 + d_2 f_1 + \frac{1}{2} d f_1 f_2 (f_1 + f_2 + e).
$$
\end{itemize}
\end{lemma}

\begin{proof}
Most of the claims are covered by \cite[Proposition
4.1]{N-gorliaison}. In any case, the proof given there shows that there is an exact sequence
$$
0 \to S (-f_1 - f_2) \to I_{C_1, X} (-f_2) \oplus I_{C_2, X} (-f_1)
\to I_{C, X} \to 0.
$$
This implies (a) and (c). Claim (b) follows from (a) by noting that the
arithmetic genus of the complete intersection $Y$ is $g_Y = \frac{1}{2} d f_1 f_2
(f_1 + f_2 + e) + 1$.
\end{proof}

\begin{corollary} \label{cor-ci-class}
The CI-liaison class of $C$ depends only on the CI-liaison  classes
of $C_1, C_2$ and the difference $f_1 - f_2$ of the degrees of the
hypersurfaces $F_1, F_2$.
\end{corollary}

\begin{proof}
Indeed, the $\mathcal N$-type resolution of $C$  involves $\cN_1
(-f_2) \oplus \cN_2 (-f_1)$, whose stable equivalence class depends
only on those of $\cN_1, \cN_2$ and the difference $f_1 - f_2$.
Hence, the claim follows from Rao's theorem \cite{rao2}, \cite[Corollary 3.14]{RTLRP} that
CI-biliaison classes are determined by the sheaf $\cN$ appearing in
the $\mathcal N$-type resolution, up to stable equivalence.
\end{proof}


\section
{Properties of Liaison Addition }

 \maketitle

We hope to use liaison addition for elucidating the structure of a
Gorenstein liaison class. To this end it is important to know under
what conditions on $C_2$ the new subscheme $C$ is in the same
G-biliaison or even G-liaison class as $C_1$.

\begin{proposition} \label{prop-rel-CI}
\begin{itemize}
\item[(a)] If $C_2$ is licci, then $C_1$ and $C$ are in the same CI-biliaison class.
\item[(b)] If $C$ and $C_1$ are in the same G-biliaison or even G-liaison class, then $C_2$ must be ACM.
\end{itemize}
\end{proposition}

\begin{proof}
(a) If $C_2$ is licci, then it has an $\mathcal N$-type  resolution
with $\cN_2$ dissoci\'e. Hence $C_1$ and $C$ have $\mathcal N$-type
resolutions with stably equivalent $\cN$.

(b) G-biliaison and even G-liaison preserve deficiency modules,  up
to twist, so the deficiency modules of $C_2$ must be all zero, i.e.\
$C_2$ is ACM.
\end{proof}

Even though the above result has a very simple proof,  it is based
on deep theorems and the links are not given explicitly (see however Section \ref{sec-licci}).

We now weaken the assumption on $C_2$.

 \begin{theorem} \label{theorem 1}
 Let $X$ be a normal arithmetically Gorenstein subscheme
 of $\mathbb P^N$ and let $C_1,C_2$ be locally Cohen-Macaulay
 codimension two subschemes of $X$.  Let $C$ be the liaison addition
 of $C_1$ and $C_2$ with respect to forms $F_1 \in I_{C_2}, F_2 \in I_{C_1}$.

 \begin{itemize}
 \item[(a)] If $C_2$ is gobilicci then $C_1$ and $C$ are in the same Gorenstein biliaison class on $X$.

 \item[(b)] If $C_2$ is glicci then $C_1$ and $C$ are in the same even Gorenstein liaison class on $X$.
 \end{itemize}
 \end{theorem}

 To prove Theorem \ref{theorem 1} we will use the fact
 that an $\mathcal N$-type resolution of $C$ is obtained
 essentially as a direct sum of $\mathcal N$-type resolutions
 of $C_1$ and $C_2$ (see Lemma \ref{lem-liai-add} (c)  above).
 Then we use criteria from the papers \cite{CH} and \cite{CDH}
 respectively characterizing subschemes in the same Gorenstein
 biliaison (resp.\ liaison) class to prove the results.

First we need an alternative form of \cite[Theorem 3.1]{CH}
characterizing Gorenstein biliaison equivalence classes.

 \begin{proposition} \label{prop 2}
Let $X$ be a normal projective arithmetically Gorenstein scheme, and
let $C_1,C_2$ be codimension two subschemes without embedded
components in $X$.  Then $C_1$ and $C_2$ are in the same G-biliaison
equivalence class if and only if they have ${\mathcal N}$-type
resolutions
 \[
 \begin{array}{c}
 0 \rightarrow {\mathcal L}_1 \rightarrow {\mathcal N}_1 \rightarrow {\mathcal I}_{C_1}(a_1) \rightarrow 0 \\
 0 \rightarrow {\mathcal L}_2 \rightarrow {\mathcal N}_2 \rightarrow {\mathcal I}_{C_2}(a_2) \rightarrow 0
 \end{array}
 \]
and there exists  an extraverti sheaf $\mathcal F$  and exact sequences (with the same
$\mathcal F$ on the left!)
 \[
 \begin{array}{ccccccccccc}
 0 & \rightarrow & {\mathcal F} & \rightarrow & {\mathcal N}_1 & \rightarrow & {\mathcal E}_1^{\vee \sigma \vee} & \rightarrow & 0 \\
 && || \\
 0 & \rightarrow & {\mathcal F} & \rightarrow & {\mathcal N}_2 & \rightarrow & {\mathcal E}_2^{\vee \sigma \vee} & \rightarrow & 0 ,
 \end{array}
 \]
where ${\mathcal E}_1$ and ${\mathcal E}_2$ are layered ACM sheaves
(see \cite{CH} for definition) of the same rank, and the rank 1
factors of the layerings of ${\mathcal E}_1$ and ${\mathcal E}_2$
are isomorphic, up to twist, in some order.  (Here $^\vee$
represents dual, and $^\sigma$ represents the syzygy sheaf, see
\cite{CDH}.)
 \end{proposition}

 \begin{proof}
This is obtained by rewriting the result of \cite[Theorem 3.1]{CH}  in
terms of the $\mathcal N$-type resolutions.  Given a sequence
 \[
 0 \rightarrow {\mathcal E} \rightarrow {\mathcal N} \rightarrow {\mathcal I}_C(a) \rightarrow 0
 \]
as in the statement of \cite[Theorem 3.1]{CH}, where $\mathcal E$ is
an ACM sheaf and $\mathcal N$ is just assumed to be coherent (note this
is neither an $\mathcal E$-type nor an $\mathcal N$-type resolution,
in spite of the notation!), we proceed as follows.

 First take a sequence
 \[
 0 \rightarrow {\mathcal L}' \rightarrow {\mathcal N}' \rightarrow {\mathcal N} \rightarrow 0
 \]
where ${\mathcal N}'$ is extraverti, which exists by \cite{HMP}, \S
2.  Then letting ${\mathcal E}'$ be the kernel of ${\mathcal N}'
\rightarrow {\mathcal I}_C(a)$, we get a new sequence
 \[
 0 \rightarrow {\mathcal E}' \rightarrow {\mathcal N}' \rightarrow {\mathcal I}_C (a) \rightarrow 0
 \]
as above, where now ${\mathcal N}'$ is extraverti.  In other words,
dropping primes, we may assume that the original $\mathcal N$ was
extraverti.

 Next take the syzygies of ${\mathcal E}^\vee$
 \[
 0 \rightarrow {\mathcal E}^{\vee \sigma} \rightarrow {\mathcal L} \rightarrow {\mathcal E}^\vee \rightarrow 0
 \]
 with $\mathcal L$ dissoci\'e.  Dualize to obtain
 \[
 0 \rightarrow {\mathcal E} \rightarrow {\mathcal L}^\vee \rightarrow {\mathcal E}^{\vee \sigma \vee} \rightarrow 0.
 \]
 Now we create a push-out diagram
 \[
 \begin{array}{ccccccccccccccccccccc}
 && 0 && 0 \\
 && \downarrow && \downarrow \\
 0 & \rightarrow & {\mathcal E} & \rightarrow & {\mathcal N} & \rightarrow & {\mathcal I}_C (a) & \rightarrow & 0 \\
 && \downarrow && \downarrow && || \\
 0 & \rightarrow & {\mathcal L}^\vee & \rightarrow & {\mathcal N}' & \rightarrow & {\mathcal I}_C (a) & \rightarrow & 0 \\
 && \downarrow && \downarrow \\
 && {\mathcal E}^{\vee \sigma \vee} & = & {\mathcal E}^{\vee \sigma \vee} \\
 && \downarrow && \downarrow \\
 && 0 && 0
 \end{array}
 \]
 The middle row is then an $\mathcal N$-type resolution of $C$.

To prove the proposition, first let $C_1, C_2$ be in the same
Gorenstein biliaison class.  Then there are sequences
 \[
 \begin{array}{cccccccccccccccccc}
 0 & \rightarrow & {\mathcal E}_1 & \rightarrow & {\mathcal N} & \rightarrow {\mathcal I}_{C_1} (a_1) & \rightarrow & 0 \\
  0 & \rightarrow & {\mathcal E}_2 & \rightarrow & {\mathcal N} & \rightarrow {\mathcal I}_{C_2} (a_2) & \rightarrow & 0
 \end{array}
 \]
as in \cite[Theorem 3.1]{CH} with the same sheaf $\mathcal N$ in the
middle.  As above, we may assume $\mathcal N$ is extraverti.  Then
performing the push-out construction for ${\mathcal E}_1$ and
${\mathcal E}_2$ as above, we get $\mathcal N$-type resolutions for
$C_1$ and $C_2$ with sheaves ${\mathcal N}_1$ and ${\mathcal N}_2$
(as the ${\mathcal N}'$ above) and exact sequences as desired with
$\mathcal F$ taken as the $\mathcal N$ above.

Conversely, given $\mathcal N$-type resolutions ${\mathcal N}_1$ and
${\mathcal N}_2$ for $C_1$ and $C_2$ respectively, and given the
sheaf $\mathcal F$ relating the two as above, for each one create a
diagram (dropping subscripts)
 \[
 \begin{array}{ccccccccccccccccccccc}
 && 0 && 0 \\
 && \downarrow && \downarrow \\
 0 & \rightarrow & {\mathcal R} & \rightarrow & {\mathcal L} & \rightarrow & {\mathcal E}^{\vee \sigma \vee} & \rightarrow & 0 \\
 && \downarrow && \downarrow && || \\
 0 & \rightarrow & {\mathcal F} & \rightarrow & {\mathcal N} & \rightarrow & {\mathcal E}^{\vee \sigma \vee} & \rightarrow & 0 \\
 && \downarrow && \downarrow \\
 && {\mathcal I}_C(a) & = & {\mathcal I}_C(a) \\
 && \downarrow && \downarrow \\
 && 0 && 0
 \end{array}
 \]
defining $\mathcal R$ as  the kernel of ${\mathcal L} \rightarrow
{\mathcal E}^{\vee \sigma \vee} \rightarrow 0$.  Then $\mathcal R$
becomes the syzygy sheaf of ${\mathcal E}^{\vee \sigma \vee}$ up to
a dissoci\'e and this is just $\mathcal E$ up to a dissoci\'e
\cite[Proposition 4.1 b]{CDH}.  So the left-hand column of these two
diagrams give the sequences of \cite[Theorem 3.1]{CH} with
${\mathcal R}, {\mathcal F}$ in place of ${\mathcal E}, {\mathcal
N}$.  Hence $C_1, C_2$ are in the same Gorenstein biliaison class.
  \end{proof}

 \begin{corollary}
 $C$ is gobilicci if and only if it has an $\mathcal N$-type resolution whose $\mathcal N$ belongs to an exact sequence
 \[
 0 \rightarrow {\mathcal E}' \rightarrow {\mathcal N} \rightarrow {\mathcal E}^{\vee \sigma \vee} \rightarrow 0
 \]
 where $\mathcal E$ and ${\mathcal E}'$ are layered ACM sheaves with the same rank 1 factors up to twist and order.
 \end{corollary}

 \begin{proof}
In the proposition we can take $C_1 = C$ and $C_2$ to be a complete
intersection.  Then ${\mathcal N}_2$ is dissoci\'e, so $\mathcal F$
becomes just ${\mathcal E}_2$ up to a dissoci\'e and we get the
desired result.

Conversely, given this sequence for $\mathcal N$, consider the
syzygy sequence for ${\mathcal E}'^\vee$,  compare its dual
 \[
 0 \rightarrow {\mathcal E}' \rightarrow {\mathcal L}^\vee \rightarrow {\mathcal E}'^{\vee \sigma \vee} \rightarrow 0,
 \]
 and apply Proposition \ref{prop 2} to see that $C$ is gobilicci.
 \end{proof}

 \bigskip

 \noindent {\em Proof of Theorem \ref{theorem 1} a.}
Let $C_1$ and $C_2$ have $\mathcal N$-type resolutions with sheaves
${\mathcal N}_1, {\mathcal N}_2$.  Assuming that $C_2$ is gobilicci,
${\mathcal N}_2$ admits a sequence
 \[
 0 \rightarrow {\mathcal E}' \rightarrow {\mathcal N}_2 \rightarrow {\mathcal E}^{\vee \sigma \vee} \rightarrow 0
 \]
 with ${\mathcal E}, {\mathcal E}'$ as above.

Then by Lemma \ref{lem-liai-add} above, $C$ has an $\mathcal N$-type
resolution with ${\mathcal N} = {\mathcal N}_1 (-f_1) \oplus
{\mathcal N}_2 (-f_2)$.  To show that $C_1$ and $C$ are in the same
Gorenstein biliaison class, we apply Proposition \ref{prop 2}.  For
simplicity, we drop the twists from the notation.

 Let
 \[
 0 \rightarrow {\mathcal E}' \rightarrow {\mathcal L} \rightarrow {\mathcal E}'^{\vee \sigma \vee} \rightarrow 0
 \]
 be the dual syzygy sequence for ${\mathcal E}'$.  Then $\mathcal L$ is dissoci\'e, so $C_1$ also has an $\mathcal N$-type resolution with ${\mathcal N}_1 \oplus {\mathcal L}$ in the middle.  We take ${\mathcal F} = {\mathcal N}_1 \oplus {\mathcal E}'$ and use the sequences
 \[
 \begin{array}{cccccccccccccccccc}
 0 \rightarrow {\mathcal F} \rightarrow {\mathcal N}_1 \oplus {\mathcal L} \rightarrow {\mathcal E}'^{\vee \sigma \vee}
 \rightarrow 0 \\
  0 \rightarrow {\mathcal F} \rightarrow {\mathcal N}_1 \oplus {\mathcal N}_2 \rightarrow {\mathcal E}^{\vee \sigma \vee}
  \rightarrow 0
 \end{array}
 \]
which show by Proposition \ref{prop 2} that $C_1$ and $C$ are in the
same Gorenstein biliaison class.   \qed

 \bigskip

 \noindent {\em Proof of Theorem \ref{theorem 1} b.}
This time we use results from \cite{CDH}.  If $C_2$ is glicci then
it has an $\mathcal N$-type resolution whose sheaf ${\mathcal N}_2
\cong {\mathcal N}_0 \oplus {\mathcal M}_0$ with ${\mathcal N}_0$
double-layered and ${\mathcal M}_0$ dissoci\'e (\cite[Corollary 5.3]{CDH}).

Let $C_1$ have an $\mathcal N$-type resolution ${\mathcal N}_1$.
Then $C$ has an $\mathcal N$-type resolution with ${\mathcal N}_1
\oplus {\mathcal N}_2$ (dropping twists).  We apply \cite[Proposition 5.1]{CDH}. Since ${\mathcal N}_0$ is double-layered, we take
the filtration given in the definition \cite[Definition 4.4]{CDH},
and insert ${\mathcal N}_1$ and ${\mathcal M}_0$ in the middle, to
satisfy the criterion of \cite[Proposition 5.1]{CDH} and show that
$C_1$ and $C$ are in the same Gorenstein liaison class.  It is an
even Gorenstein liaison class because the sheaf in the middle of the
filtration is ${\mathcal N}_1$ and not ${\mathcal N}_1^{\sigma
\vee}$.  \qed


 \section{Licci subschemes} \label{sec-licci}

It is interesting to see that  the
links  needed in  Proposition \ref{prop-rel-CI} can be described in a very concrete way.  We begin with the
following preliminary tool.

\begin{proposition} \label{la is liaison}
Let $X$ be a normal arithmetically Gorenstein subscheme  of $\mathbb
P^n$ and let $C_1$ and $C_2$ be  codimension two subschemes of $X$.
Choose $B,F \in I_{C_1}$ and $A,G \in I_{C_2}$ such that $AF$ and
$BG$ form a regular sequence.  We make the following codimension two
subschemes:

\begin{enumerate}
\item $C_1'$ is the residual to $C_1$ in the complete intersection $(F,B)$;

\item $C_2'$ is the residual to $C_2$ in the complete intersection $(G,A)$;

\item $C$ is the liaison addition subscheme defined by $I_C = G \cdot I_{C_1} + F \cdot I_{C_2}$;

\item $C'$ is the liaison addition subscheme defined by $I_{C'} = A \cdot I_{C_1'} + B \cdot I_{C_2'}$.
\end{enumerate}
Then $C$ and $C'$ are directly linked by the complete intersection $(AF,BG)$.

\end{proposition}

\begin{proof}
By hypothesis we have
\[
\begin{array}{rcl}
(F,B) : I_{C_1} & = & I_{C_1'} \\
(G,A) : I_{C_2} & = & I_{C_2'}
\end{array}
\]
We have to show that
\begin{equation} \label{to show}
(AF,BG) : (G\cdot I_{C_1} + F \cdot I_{C_2}) = A \cdot I_{C_1'} + B \cdot I_{C_2'}.
\end{equation}
We proceed in two steps.  First we show the inclusion $\supseteq$.
To this end, let $H \in I_{C_1'}$.  By construction we have $H \cdot
I_{C_1} \subset (F,B)$.  Then it follows that
\[
AH \cdot (G\cdot I_{C_1} + F \cdot I_{C_2}) \subset (AF,BG)
\]
since
\[
AHG \cdot I_{C_1} \subset (AGF, AGB) \subset (AF,BG)
\]
and clearly
\[
AHF \cdot I_{C_2} \subset (AF,BG)
\]
Now let $K \in I_{C_2'}$.  In a completely analogous way we have
\[
BK \cdot (G\cdot I_{C_1} + F \cdot I_{C_2}) \subset (AF,BG).
\]
We thus have shown the inclusion $\supseteq$ of (\ref{to show}).

Now, all ideals under consideration are the saturated ideals of
codimension  two subschemes of $X$.  Hence the equality of (\ref{to
show}) will be established if we can show that both sides define
schemes of the same degree.  Let $f = \deg F$, $g = \deg G$, $a =
\deg A$,  $b = \deg B$, and $d = \deg X$.  Then
\[
\begin{array}{rcl}
\deg C & = & \deg C_1 + \deg C_2 + d fg \\
\deg C' & = & \deg C_1' + \deg C_2' + d  ab
\end{array}
\]
and we know $\deg C_1 + \deg C_1' = d bf$ and $\deg C_2 + \deg C_2' = d ag$, so
\[
\begin{array}{rcl}
\deg C + \deg C' & = & d (bf + ag + fg + ab) \\
& = & d (a+f) (b+g),
\end{array}
\]
which completes the proof.
\end{proof}

In the following corollary, the fact that $C$ is evenly CI-linked to
$C_1$  follows immediately from Rao's theorem, as noted above, since
the $\mathcal N$-type resolutions give bundles that are stably
equivalent.  The content of this corollary is that we can follow the
links in such a precise way.

 \begin{corollary} \label{licci la}
Let $X$ be a normal  arithmetically Gorenstein subscheme of $\mathbb P^n$ and let $C_1, C_2$ be codimension two subschemes of $X$.  Assume that $C_2$ is {\em licci}, with $r$ minimal generators.  Choose $F \in I_{C_1}$ and $G \in I_{C_2}$ such that $F$ and $G$ form a regular sequence, and let $C$ be the liaison addition subscheme defined by the saturated ideal $I_C = G \cdot I_{C_1} + F \cdot I_{C_2}$.  Then $C$ is  CI-linked to $C_1$ in an even number of steps.  Furthermore,  there is a sequence of subschemes
\[
C, Z_{r-1}, Y_{r-1}, Z_{r-2}, Y_{r-2}, \dots, Z_2, Y_2, Z_1, Y_1, D, C_1
\]
where
\begin{itemize}
\item[(a)] any two consecutive subschemes in the sequence are directly linked (by complete intersections that we will specify);
\item[(b)] for $i \geq 2$, each $Y_i$ is obtained as the liaison addition of $C_1$ with a licci subscheme with $i$ minimal generators;
\item[(c)] Each $Z_i$ is obtained as the liaison addition of a (fixed) subscheme directly linked to $C_1$ with a licci subscheme;
\item[(d)] $Y_1$ is a {\em basic double link} of $C_1$, i.e. the liaison addition of $C_1$ with the trivial subscheme.
\end{itemize}
 \end{corollary}

 \begin{proof}
 By Rao's theorem, a codimension two licci subscheme $Y$ of $X$ has a minimal free $R/I_X$-resolution of the form
 \[
 0 \rightarrow L_2 \rightarrow L_1 \rightarrow I_{Y} \rightarrow 0.
 \]
By a standard trick due to
Gaeta (in modern language this is shown via mapping cones -- cf.\
\cite{book}),  we have the following possibilities for linking $Y$.
Suppose that $F_1,F_2 \in I_{Y}$ are a regular sequence, linking $Y$
to a residual subscheme $Y'$ of $X$.

 \begin{enumerate}
\item If $F_1$ and $F_2$ are both minimal  generators
of $I_{Y}$ then neither is a minimal generator of $I_{Y'}$.   In
this case $I_{Y'}$ has one less minimal generator than does $I_{Y}$.

 \item If $F_1$ is a minimal generator of $I_{Y}$ but $F_2$ is
 not, then $F_1$ is again a minimal generator of $I_{Y'}$, but $F_2$ is not.
 In this case $I_{Y'}$ has the same number of minimal generators as does $I_{Y}$.

 \item If neither $F_1$ nor $F_2$ are minimal generators of $I_{Y}$
 then both are minimal generators of $I_{Y'}$.  In this
 case $I_{Y'}$ has one more minimal generator than does $I_{Y}$.

 \end{enumerate}

Now, in our situation, we will show that $C$ can be linked in two
steps to a codimension two subscheme $C''$ that is the liaison
addition of $C_1$ and a licci subscheme $C_2''$ whose ideal has one
fewer minimal generator than does  $I_{C_2}$.  The result will then
follow by induction.

Choose $B \in I_{C_1}$ and $A \in I_{C_2}$ such that $A$  and $B$
form a regular sequence, and such that furthermore $A$ is a minimal
generator of $I_{C_2}$.  As before, let $C_1$ be directly linked to
$C_1'$ by $(F,B)$ and let $C_2$ be directly linked to $C_2'$ by
$(G,A)$.  Then by Proposition \ref{la is liaison}, $C$ is directly
linked via $(AF,BG)$ to the liaison addition subscheme $C'$
corresponding to the ideal $A \cdot I_{C_1'} + B \cdot I_{C_2'}$.

Now replace the data $(C_1,C_2,F,G,B,A)$ by the data
$(C_1',C_2',B,A,F,A')$ where $A'$ is a minimal generator of
$I_{C_2'}$.  We get that ${C_2'}$ is directly linked by $(A,A')$ to
a subscheme $C_2''$, $C_1'$ is directly linked by $(B,F)$ back to
$C_1$, and $C'$ is directly linked by $(A'B, FA)$ to a subscheme
$C''$ defined by the liaison addition $I_{C''} = A' \cdot I_{C_1} +
F \cdot I_{C_2''}$.


Now we consider the number  of minimal generators of $I_{C_2}$ and
$I_{C_2''}$.  Since $G$ may or may not have been a minimal generator
of $I_{C_2}$, while $A$ was a minimal generator, we have two
possibilities.

If $G$ was a minimal generator  of $I_{C_2}$ then $I_{C_2'}$ has one
fewer minimal generator than does $I_{C_2}$, but then $I_{C_2''}$
has the same number of minimal generators as $I_{C_2'}$, which is
one less than $I_{C_2}$.

If $G$ was not a minimal generator  of $I_{C_2}$ then $I_{C_2'}$ has
the same number of minimal generators as $I_{C_2}$, but then $A$ and
$A'$ are both minimal generators of $I_{C_2'}$, so $I_{C_2''}$ has
one fewer minimal generator.  Note that neither $A$ nor $A'$ are
minimal generators of $I_{C_2''}$.

By induction, we arrive in an even  number of steps to the liaison
addition of $C_1$ and a complete intersection, $C_2$.  As we have
seen above, but using the notation of Proposition \ref{la is
liaison},  the polynomial $G \in I_{C_2}$ used in the liaison
addition is not one of the minimal generators of $I_{C_2}$.  One
more link as we did above results in the subscheme $C'$ consisting
of the liaison addition of $C_1'$ with another complete
intersection, but this time the polynomial $A \in I_{C_2'}$ used in
the liaison addition {\em is} a minimal generator of $I_{C_2'}$.  We
will write $I_{C_2'} = (A,A')$, where $A$ and $A'$ have no common
factor.

So at this stage we are considering  the subscheme $C'$ which is a
liaison addition of the form
\[
A \cdot I_{C_1'} + B \cdot (A',A).
\]
This is linked by the complete  intersection $(A'B,FA)$ to the
subscheme $C''$ defined by $A' \cdot I_{C_1} + (F)$, since $C_2'$ is
linked by $(A,A')$ to the trivial ideal $R/I_X$.  But $C'$ is
precisely a basic double link ideal, which is linked in two steps to
$C_1$.
 \end{proof}

 \begin{remark}
In the proof of Corollary \ref{licci la}, we
 used in a heavy way the theory of liaison for codimension two licci ideals.
 We wonder about the following questions.

 \begin{enumerate}
 \item Since Proposition \ref{la is liaison} did not
 assume that $C_2$ was licci, it should have additional applications.
 If we start with {\em any} liaison addition of subschemes $C_1$ and $C_2$,
 can we explicitly link it in a finite number of steps to a suitable liaison addition
 of a {\em minimal} element in the even liaison class of $C_1$ and a
 minimal element in the even liaison class of $C_2$?

 \item Liaison addition and basic double linkage have been
 developed for higher codimension  \cite{BM4}, \cite{GM4}, \cite{N-gorliaison} in a
 way that is very similar to the codimension two picture.  Can
 Proposition \ref{la is liaison} and Corollary \ref{licci la} also be extended to higher codimension?
 \end{enumerate}

 \end{remark}


\section{Curves on quadric threefolds in $\PP^4$ -- toward a Lazarsfeld-Rao-type structure}

There is a beautiful structure theorem for  an even liaison class of
codimension two subschemes.  It was discovered for curves in
$\mathbb P^3$ by Martin-Deschamps and Perrin \cite{MP} and for
codimension two subschemes in $\mathbb P^n$ by Ballico, Bolondi and
Migliore \cite{BBM}, based on a conjecture of Harris and a special case proved by Lazarsfeld and Rao.
 It has been extended to codimension two
subschemes of arithmetically Gorenstein varieties in
\cite{Bolondi-Migliore-manuscripta} and in a more general way in
\cite{N-gorliaison} and in \cite{RTLRP}, but always for even
CI-liaison.  It was pointed out in \cite{book} that extending this
property to Gorenstein liaison will be difficult.  Here, we study the question in two special cases.
\medskip

\noindent
{\bf Nonsingular quadrics}
\smallskip

Let $X$ be a non-singular quadric 3-fold in ${\mathbb P}^4$.   The
only surfaces on $X$ are complete intersections, so for curves in
$X$, Gorenstein biliaison is equivalent to CI-biliaison, which is
also equivalent to even CI-liaison.  The even CI-liaison class of a
curve $C$ is determined by a triple $(M,P,\alpha)$, where $M =
H^1_*({\mathcal I}_C)$ is the Rao module of $C$, $P$ is a maximal
Cohen-Macaulay module over the homogeneous coordinate ring of $X$,
say $S = H^0_* ({\mathcal O}_X) = R/Q$ (where $Q$ is the defining
polynomial of $X$), and $\alpha : P^\vee \rightarrow M^* \rightarrow
0$ is a surjective map of graded $S$-modules (\cite[Corollary 4.3]{RTLRP}).  This triple is determined up to isomorphism and
shift of degrees for $M$, up to stable equivalence and (the same)
shift of degrees for $P$, and compatible maps $\alpha$.

On the other hand, the even G-liaison class of a curve on $X$ is
determined by the Rao module alone (up to shift) -- cf.\ \cite[Theorem 6.2]{CDH}.

For each CI-biliaison equivalence class we have the traditional
Lazarsfeld-Rao theorem, \cite[Theorem 3.4]{RTLRP}.  Each even
Gorenstein liaison class is a union of CI-biliaison classes, so we
can ask what kind of structure the even Gorenstein liaison class can
have.  That is, how are the CI subclasses related?

We will show that for ACM curves, which form one even Gorenstein
liaison class since their Rao modules are zero, we can obtain all the non-licci curves  by a combination of CI-biliaisons and liaison additions with a
line, starting from a line.   On the other
hand,  for an even Gorenstein liaison class of curves with
 non-zero Rao module, it is not possible to obtain them
all from a single one, or even from minimal ones,  by  CI-biliaisons and
liaison additions with ACM curves (see  Remark
\ref{rem-caution} and Example \ref{ex-counter}).

\begin{theorem} \label{thm-sm-quad}
Every CI-biliaison class of non-licci ACM curves on the non-singular quadric 3-fold $X$ contains a minimal curve $C$ that can be obtained by liaison addition with a line from a minimal curve of lower degree in another such class, unless $C$ is already a line.
\end{theorem}

The proof requires some preparation.
Let $L \subset X$ be a line. Let $\cE_0$ be the locally free sheaf
defined by the minimal $\mathcal N$-type resolution of $L$:
$$
0 \to \cO_X \to \cE_0 \to \cI_{L, X} (1) \to 0.
$$
Then, according to \cite{CDH},  each non-licci CI-liaison class of ACM curves on $X$ corresponds via its
$\mathcal N$-type resolution to the stable equivalence class of one of the
sheaves
$$
\cN_{0, a_2,\ldots,a_r} := \cE_0 \oplus \cE_0 (-a_2) \oplus \cE_0
(-a_2 - a_3) \oplus \ldots \oplus \cE_0 (-a_2 - \ldots - a_r)
$$
where $a_2,\ldots,a_r \geq 0$ are integers. Note that $\cN_0 = \cE_0$.

Using liaison addition we first construct curves that we will then show to be minimal in their CI-biliaison classes.

\begin{lemma} \label{lem-res}
For each $a_2,\ldots,a_r \geq 0$ there is a curve $C_{0, a_2,\ldots,a_r}$ with $\mathcal N$-type resolution
\begin{eqnarray*}
\lefteqn{
0 \to \left ( \cO_X \oplus \cO_X^2 (-a_2) \oplus \ldots \oplus \cO_X^2
(-a_2 - \ldots - a_r) \right  ) (-r) \to } \\
 & & \hspace*{6cm}  \cN_{0, a_2,\ldots,a_r}(-r) \to
\cI_{C_{0,a_2,\ldots,a_r}, X} \to 0
\end{eqnarray*} and with $\mathcal E$-type resolution:
\begin{eqnarray*}
\lefteqn{
 0 \to \cN_{0, a_2,\ldots,a_r}(-r-1) \to } \\
 & & \left (\cO_X^3 \oplus \cO_X^2(-a_2) \oplus \ldots \oplus
 \cO_X^2(-a_2 - \ldots - a_r) \right ) (-r) \to
 \cI_{C_{0,a_2,\ldots,a_r}, X} \to 0.
\end{eqnarray*}

Furthermore, for $r \geq 2$ the curve $C_{0, a_2,\ldots,a_r}$ can be obtained by liaison addition with a line from $C_{0, a_3,\ldots,a_r}$.
\end{lemma}

\begin{proof}
We use induction on $r$. For $r=1$, we take
$C_0$ to be a line $L$. Its $\cN$-  and  $\mathcal E$-type resolution are well-known.
If $r \geq 2$, then the $\mathcal E$-type resolution of $C_{0, a_3,\ldots,a_r}$
shows that its homogenous ideal
contains an element $F_2$ of degree $r-1$. Let $F_1 \in I_{L, X}$ be
an element of degree $a_2 + 1 \geq 1$, such that $\{F_1, F_2\}$ is a
regular sequence. Then Lemma \ref{lem-liai-add} shows that the curve $C_{0, a_2,\ldots,a_r}$ defined by the liaison addition ideal $F_2 \cdot
I_L + F_1 \cdot I_{C_{0, a_3,\ldots,a_r}}$ has the required $\mathcal N$-type resolution.
The  $\mathcal E$-type
resolution of the curve $C_{0, a_2,\ldots,a_r}$ can be obtained using the syzygy sequence
\[
0 \to \cE_0 (-1) \to \cO_X^4 \to \cE_0 \to 0
\]
and the method of converting an $\cN$-type resolution to an $\cE$-type resolution described in \cite[Proposition 4.3(a)]{CDH}.

This proves the lemma, since all these curves have been constructed by liaison addition with a line, starting with $C_0 = L$ which is a line.
\end{proof}

We now show that the constructed curves are minimal in their CI-biliaison classes.

\begin{lemma} \label{lem-min}
Each curve $C_{0, a_2,\ldots,a_r}$ described in Lemma \ref{lem-res} is minimal in its CI-biliaison class.
\end{lemma}

\begin{proof}
For this lemma we will change the above notation and rewrite the $\cN$-type resolution of $C = C_{0, a_2,\ldots,a_r}$ as
$$
0 \to \cO_X^{2 r_1 -1} \oplus \cO_X^{2 r_2} (-b_2) \oplus \ldots \oplus \cO_X(- b_m)^{2 r_m} \to \cN \to \cI_{C, X}(r) \to 0
$$
where $\cN = \cE_0^{r_1} \oplus \cE_0(-b_2)^{r_2}  \oplus \ldots \oplus \cE_0(-b_m)^{r_m}$ is of total rank $r = 2 \sum r_i$ and $0 = b_1 < b_2 < \ldots < b_m$. To prove that $C$ is minimal, it is enough to show that if
$$
0 \to \bigoplus_{i=1}^{2r-1} \cO_X (-c_i) \to \cN \to \cI_{D, X}(s) \to 0
$$
is the $\cN$-type resolution of any other curve $D$ on $X$, then
$$
\sum_{i=1}^{2r-1} c_i \geq 2 \sum_{j=1}^m r_j b_j.
$$

We may assume that $c_1 \leq \cdots \leq c_{2r-1}$.
It follows from the existence of $D$ that for any subset $J$ of $\{1,\ldots,2r-1\}$, the cokernel of the map
$$
\alpha_J: \bigoplus_{i \in J} \cO_X (-c_i) \to \cN
$$
is torsion free. Let $J := \{c_i \s c_i < b_2\}$. Then the image of $\alpha_J$ lands inside the sheaf $\cE_0^{r_1}$, so the rank of $\alpha_J$ must be less than $2 r_1$. This means that for all $i \geq 2 r_1$, we have $c_i \geq b_2$.

Next,  let $J := \{c_i \s c_i < b_3\}$. Then the image of $\alpha_J$ lands inside the sheaf $\cE_0^{r_1} \oplus \cE_0(- b_2)^{r_2}$, and the same argument shows that, for all $i \geq 2 (r_1 + r_2)$, we must have $c_i \geq b_3$.

Continuing in this fashion, the inequalities on the $c_i$'s and $b_j$'s show that $\sum c_i \geq 2 \sum r_j b_j$, as required.
\end{proof}

The above theorem follows now easily.

\begin{proof}[Proof of Theorem \ref{thm-sm-quad}]
Each non-licci CI-biliaison class of ACM curves on $X$ corresponds to one of the sheaves $\cN_{0, a_2,\ldots,a_r}$ of rank $\geq 2$. Thus, the theorem follows by combining Lemmas \ref{lem-res} and   \ref{lem-min} and the Lazarsfeld-Rao property for CI-biliaison classes.
\end{proof}

\begin{remark} \label{rem-caution}
In the case  of the traditional
Lazarsfeld-Rao property, Bolondi and Migliore
\cite[Corollary 4.10]{Bolondi-Migliore-manuscripta} have shown that
one can get from a minimal scheme in a codimension two CI-liaison
class to any other (up to flat deformation) by liaison addition with
a  licci scheme. In general, the analogous result is not true in even Gorenstein
liaison classes.
\end{remark}

\begin{example} \label{ex-counter}
To see that a method similar to that used for ACM curves cannot
work with non-ACM classes of curves on $X$, consider the even
Gorenstein liaison class of curves with Rao module $M = k$.  Among
these there is one CI-biliaison class corresponding to a triple
$(M,P,\alpha)$, where $P$ is a non-free maximal Cohen-Macaulay
module over $S$, and $ \alpha : P^\vee \rightarrow M^* \rightarrow
0$.  There is another CI-biliaison class corresponding to the triple
$(M,S,\beta)$, where we think of $S$ as the free rank 1 $S$-module
and $\beta : S \rightarrow k \rightarrow 0$ the natural map \cite[Corollary 4.3]{RTLRP}.  Since
liaison addition acts as direct sum on Rao modules and on $\cN$-type resolutions, one sees easily that it also acts as direct
sums on the triples $(M,P,\alpha)$.  An ACM curve will have triple
$(0,Q,\sigma)$.  Adding this to the first curve above will give a
triple $(M,P\oplus Q,\alpha)$, where $\alpha$ acts by 0 on the $Q$
factor.  Similarly, adding to the second will give $(M,S \oplus Q,
\beta)$, with $\beta$ acting as 0 on the $Q$-factor.  So it is clear
that no combination of liaison additions with ACM curves will ever
connect these two types of curves with Rao module $k$, because in
one case $k$ is covered by a non-free maximal Cohen-Macaulay module
and in the other case by a free maximal Cohen-Macaulay module.
\end{example}


\noindent
{\bf Singular quadrics with one double point}
\smallskip

Let $X$ be the singular quadric threefold in $\mathbb P^4$ having
just one double point.  In this case two curves are in the same
Gorenstein biliaison equivalence class if and only if their Rao
modules are isomorphic up to shift (\cite[Theorem 6.2]{CH}).  It
follows that Gorenstein biliaison is the same as even Gorenstein
liaison in this case.  Indeed, we know in general that any
Gorenstein biliaison is an even Gorenstein liaison.  Conversely, if
two curves are in the same even Gorenstein liaison class, then their
Rao modules are isomorphic up to shift, and so by the theorem above
they are equivalent for Gorenstein biliaison.

Thus, having the operation of Gorenstein biliaison available, one
might hope, as in the case of the traditional Lazarsfeld-Rao
theorem, that in any Gorenstein biliaison equivalence class, every
curve could be obtained by a  finite sequence of ascending
Gorenstein biliaisons from one, or a small number of ``minimal"
curves in the class. The following example shows that this is too ambitious.

 \begin{example} \label{ex1}
For ACM curves on $X$, the natural choice
for minimal curves would be a line or a conic.  There are three
types of lines: those contained in a $D$-plane, those contained in
an $E$-plane, or those passing through the double point of $X$,
where $D$ and $E$ refer to the rulings over a general hyperplane
section $Q$ of $X$ (which is a smooth quadric surface), and we think
of $X$ as a cone over $Q$.  These lines have $\mathcal N$-type
resolution using ${\mathcal E}_1, {\mathcal E}_1'$ or ${\mathcal
I}_D \oplus {\mathcal I}_E$ respectively, in the notation of
\cite[Theorem 6.2]{CH}.  However, we will exhibit here an infinite
sequence of ACM  curves that do not admit
any descending Gorenstein biliaisons.

Take two $D$-planes in $X$, say $D_1$ and $D_2$.  They meet only at
the singular point $P$ of $X$.  Take curves $C_1 \subseteq D_1$,
$C_2 \subseteq D_2$ of degrees $d$ and $e$, respectively, each
passing simply through the common point $P$ of the two planes.
Consider  the exact sequence
 \[
 0 \rightarrow I_{C_1} \cap I_{C_2} \rightarrow I_{C_1} \oplus I_{C_2} \rightarrow I_{C_1} + I_{C_2} \rightarrow 0.
 \]
Note  that $I_{C_1} + I_{C_2} = I_P$, and that $C_1$ and $C_2$ are
both ACM, being plane curves.  Then
sheafifying this sequence and taking cohomology, it is immediate to
see that $C = C_1 \cup C_2$ is an ACM
curve of degree $d+e$.  Assuming that $e \geq 2$ and $d \geq e+2$,
we will show that $C$ does not admit any descending Gorenstein
biliaison in $X$.

So suppose that $C$ admits a descending Gorenstein biliaison on an
ACM surface $Y$ in $X$.  We distinguish
four cases.

 \medskip

 \begin{itemize}
 \item[\underline{Case 1}:]
$Y = D_1 \cup D_2 \cup F$ for some other surface $F$.  This is
impossible, because there is nothing in $F$ to subtract a hyperplane
section from.

 \medskip

 \item[\underline{Case 2}:]
$Y = D_1 \cup F$, where $F$ is some other surface containing $C_2$,
but not containing $D_2$.  Letting $(a,b)$ be the bidegree of $F$ on
$X$, i.e. $F \sim aD + bE$ on $X$, we see that $b \geq e$ by
intersecting $F$ with $D_2$.  On the other hand, since $C_2$ passes
only simply through $P$, we must have $a>0$.  Then the degree of $F$
is $a+b > e$, and we cannot subtract a hyperplane section of $F$
from $C_2$.

 \medskip

 \item[\underline{Case 3}:]
 $Y = F \cup D_2$.  This is similar to Case 2.

 \medskip

 \item[\underline{Case 4}:]
$Y$ does not contain either $D_1$ or $D_2$.  Then $Y$ has bidegree
$(a,b)$ with $b \geq d$, and since $Y$ is an ACM surface on $X$,
we must have $|a-b| \leq 1$.    Hence
 \[
 \deg Y = a+b \geq (b-1)+b \geq (d-1)+d > e+d,
 \]
because of our hypotheses on $d$ and $e$.  Hence we cannot bilink
down on $Y$.

 \end{itemize}

 \end{example}

 \bigskip

 \begin{remark}
As in the case of the non-singular quadric threefold, one could ask
whether every ACM curve can be obtained by
a succession of complete intersection bilinks and liaison additions,
for example with a line.  On the non-singular quadric threefold $X$,
there is only one non-trivial indecomposable ACM sheaf, up to twist,
namely ${\mathcal E}_0$, and it corresponds to a line.  Thus any ACM
sheaf can be obtained by adding direct sums of twists of this one to
a dissoci\'e sheaf and this explains why the method shown above
works in this case.

On the singular quadric threefold, there are two infinite sequence
of indecomposable ACM sheaves, ${\mathcal E}_\ell$ and ${\mathcal
E}_\ell'$ for $\ell = 1,2,\dots$ \cite[proof of Theorem 6.2]{CH}, so in order to formulate an
analogous result, one would have to allow (at least) liaison
additions with plane curves of all degrees in both $D$-planes and
$E$-planes, either passing or not passing through the singular point
$P$.

 \end{remark}

 \bigskip

 \begin{example}
Now we consider curves with Rao module $k$.  In $\mathbb P^4$ the
curves with minimal leftward shift of $k$, namely $k$ in degree 0, have
been classified \cite{lesperance}, \cite[Proposition 4.1]{SEG}.  They
exist in any degree $d \geq 2$, and the general such curve is the
disjoint union of a plane curve of degree $d-1$ and a line not
meeting the plane of the first curve.  Curves of this kind for every
$d \geq 2$ can be found on the singular quadric threefold $X$, so we
take these as the minimal curves.  One might hope that curves whose
Rao module has a shift into positive degrees of the module $k$ could
be obtained by ascending Gorenstein biliaison from these minimal
curves.  We show this is not the case by exhibiting some non-singular curves of degree 5 and genus 0 on $X$ that have Rao module $k$ in degree 1 and do not admit any descending Gorenstein biliaison on $X$.

We begin by recalling some basic facts about degree 5 and genus 0 curves $C$ in $\PP^4$ \cite[Example 4.3]{SEG}. Such curves can be obtained by generic projection from the rational normal curve in $\PP^5$. As long as the curve is non-degenerate (i.e.\ not contained in any $\PP^3$), we find from Riemann-Roch that $h^0(\cI_C(2)) = 4$. Taking two general quadric hypersurfaces containing $C$, the curve $C$ will be contained in their intersection $Y$, a degree 4 Del Pezzo surface, which will be non-singular provided $C$ is general. Conversely, on the Del Pezzo surface $Y$, we can find  smooth curves in the divisor classes $(2; 1, 0^4), (3; 2, 1^2, 0^2), (4; 2^3, 1, 0)$, and $(5; 3, 2^3, 1)$ (and their permutations), which we will denote by $C_1, C_2, C_3, C_4$, respectively. Here we use the standard notation for divisor classes on $Y$ \cite[Notation 3.3]{SEG}. There is no difference between these curves as curves in $\PP^4$. However, as curves on $Y$ they are distinguished by thei
 r divisor classes on $Y$.

Next we observe that each general non-singular, non-degenerate curve $C$ in $\PP^4$ of degree 5 and genus 0 has a unique trisecant. Indeed, any trisecant of $C$  must lie in every quadric hypersurface containing $C$ and therefore on the Del Pezzo surface $Y$. Then, checking each of the lines on $Y$, one finds exactly one trisecant for each $C_i$ on $Y$.

We can also find  non-singular degree 5 genus 0 curves $C_5$ on the rational cubic scroll $S$ in the divisor class $(4; 3)$, and this is the only possibility. This curve $C_5$ meets a fiber $F = (1; 1)$ of the ruling in one point, and it meets the exceptional curve $E = (0; 1)$ three times. Thus, the ruling determines an isomorphism of $C_5$ to the projective line $E$. This isomorphism is uniquely determined by the three intersection points of $C_5$ with $E$ (which are necessarily distinct), and we can recover the surface $S$ as the closure of the union of lines joining corresponding points on $C_5$ and $E$. This shows that $C_5$ lies on a unique cubic scroll. Since we saw earlier that each general degree 5 genus 0 curve in $\PP^4$ has a unique trisecant, it follows that every such curve $C$ is contained in this way in a cubic scroll.

On the cubic scroll, the divisor class $C - H = (2; 2)$ contains a union of two fibers of the rulings. The passage from $C-H$ to $C$ is a Gorenstein biliaison in $\PP^4$, hence $C$ has Rao module $k$ in degree 1, and as a curve in $\PP^4$ it does admit a descending Gorenstein biliaison \cite[Example 4.3]{SEG}.

Now we consider non-singular degree 5 genus 0 curves $C$ on the singular quadric threefold $X$ in $\PP^4$. There are cubic scrolls $S$ in $X$, having bidegrees $(2, 1)$ and $(1, 2)$. If $C$ is in $S$, then its projection $\pi (C)$ from the singular point $P$ of $X$ to a general hyperplane section $Q$, which is a non-singular quadric surface in $\PP^3$, is isomorphic to $C$ because the projection maps $S$ to $Q$ birationally, blowing up the point $P$ and blowing down the ruling through $P$, which meets $C$ just once. Hence $\pi (C)$ is a curve of bidegree $(1, 4)$ or $(4, 1)$ on $Q$.

There are also Del Pezzo surfaces $Y$ on $X$, having bidegree $(2, 2)$. The intersections of the two families of planes on $X$ with $Y$ are conics adding to a hyperplane section of $Y$, and without loss of generality we can take these to be $\Gamma = (1; 1, 0^4)$ and $\Gamma' = (2; 0, 1^4)$. For any curve $C$ on $X$, the two intersection numbers $C.\Gamma, C.\Gamma'$ will give the bidegree of the projection $\pi (C)$ of $C$ onto the hyperplane section $Q$. Thus we find that $\pi (C_1)$ and $\pi (C_2)$ have bidegree $(4, 1)$ or $(1, 4)$, whereas $\pi (C_3)$ and $\pi (C_4)$ have bidegree $(3, 2)$ or $(2, 3)$. Since we saw above that a  $(5,0)$ curve on a cubic scroll must project onto a curve with bidegree $(1, 4)$ or $(4, 1)$ on $Q$, the curves $C_3$ and $C_4$ cannot be contained in any cubic scroll on $X$.

Now, finally, we show that neither $C_3$ nor $C_4$ admit any descending biliaison on $X$. Indeed, let $C$ be one of these two curves and suppose that $C$ is contained in an ACM surface $T$ in $X$ and that $C - H$ is effective on $T$. Then $\deg (C - H) = 5 - \deg T$ must be at least 2, since any curve of degree 1 is ACM. So $\deg T \leq 3$. The degree cannot be 2, since $C$ is not contained in a hyperplane. We conclude that $\deg T = 3$. However, the only irreducible surfaces of degree 3 in $X$ are the cubic scroll, which does not contain $C$, and the cone over a twisted cubic, which contains no non-singular curves of degree 5 and genus 0. Hence a descending Gorenstein biliaison of $C$ is not possible on $X$.

Note of course that either curve $C_3$ or $C_4$ is contained in a unique cubic scroll in $\PP^4$, but this argument shows that in this case the cubic scroll lies outside of $X$, and intersects $X$ only in the curve $C_i$ ($i = 3, 4$) together with its trisecant.
\end{example}

We close this section by wondering if the concept of liaison addition can be extended. This could  potentially be useful to address some of the problems we encountered above.

\begin{remark} \label{rem-liai-ext}
Let $C_1, C_2$ be two codimension 2 subschemes of some projective scheme $X \subset \PP^n$ with Rao modules $M_1$ and $M_2$. Then the liaison addition with respect to hypersurfaces of degrees $d_1, d_2$ is a curve with Rao module $M_1 (-d_2) \oplus M_2 (-d_1)$. Assume now that $N$ is any graded module corresponding to an extension
$$
0 \to M_2 (-d_1) \to N \to M_1 (-d_2) \to 0.
$$
Is it then possible to construct a curve $C$ starting directly from $C_1$ and $C_2$ such that the Rao module of $C$ is isomorphic to $N$? If the answer is affirmative, this could possibly provide a natural extension of liaison addition and it would be justified to call the curve $C$  a {\em liaison extension} of $C_1$ and $C_2$.
\end{remark}

\section{Conclusion}

Our motivation for the work in this paper was to investigate Questions 1 and 2 from the introduction: What is the structure of a Gorenstein biliaison class or an even Gorenstein liaison class of codimension 2 subschemes of an arithmetically Gorenstein subscheme? To address this question we employ the idea of liaison addition. Liaison addition has been used to investigate CI-liaison classes. Here we show that it can also be used to study Gorenstein liaison classes.

Our two main test cases have been the non-singular quadric 3-fold and the singular quadric 3-fold with one double point in $\PP^4$, because in each case we have available explicit descriptions of CI-biliaison classes as well as of Gorenstein biliaison and even liaison classes of curves.

We found a satisfactory answer for ACM curves on the non-singular quadric 3-fold where the non-licci ACM curves can all be obtained starting from a single line by successive liaison additions with a line and CI-biliaisons. One could perhaps get an analogous result for ACM curves on the singular quadric 3-fold, using liaison additions with plane curves of all degrees in both families. This gives some hope for liaison addition as a key operation in explaining the structure of an even Gorenstein liaison class. However, simple examples show that this alone is not sufficient to deal with the case of non-ACM curves. Thus, we wonder if liaison addition can be extended to liaison extension.

On the singular quadric 3-fold we have available the method of Gorenstein biliaison, and one might hope to reach any curve by ascending Gorenstein biliaison from a suitable class of ``minimal'' curves. For curves with non-zero Rao module, the natural definition would be those whose Rao module has the left-most shift. But examples show that with this definition there are non-minimal curves with no descending biliaisons. For ACM curves, we have found infinitely many classes of curves with no descending biliaisons, so there does not appear to be a suitable class of ``minimal'' curves from which all others can be obtained by ascending biliaisons.

In summary, some new ideas will be necessary to give satisfying answers to Questions 1 and 2 in general.

\end{document}